\documentclass[12pt]{amsart}
\usepackage{amsmath}
\usepackage{latexsym}
\usepackage{amscd}
\DeclareMathOperator{\im}{im}
\newcommand{\bd}{\partial}
\title[Splitting homomorphisms and the Geometrization Conjecture]
{Splitting homomorphisms and the \\ Geometrization Conjecture}
\author{Robert Myers} 
\address{Department of Mathematics, Oklahoma State University, Stillwater, OK 74078}
\email{myersr@math.okstate.edu}
\subjclass{Primary: 57N10; Secondary: 57M50, 57M10, 57M40, 57N12}
\keywords{3-manifold, Heegaard splitting, splitting homomorphism, Geometrization 
Conjecture} 

\newtheorem{thm}{Theorem}[section]
\newtheorem{prop}[thm]{Proposition}
\newtheorem{cor}[thm]{Corollary}
\newtheorem{lem}[thm]{Lemma}

\begin{document}

\begin{abstract} This paper gives an algebraic conjecture which is shown to be 
equivalent to Thurston's Geometrization Conjecture for closed, orientable 
3-manifolds. It generalizes the Stallings-Jaco theorem which established a similar 
result for the Poincar\'{e} Conjecture. The paper also gives two other algebraic 
conjectures; one is equivalent to the finite fundamental group case of 
the Geometrization Conjecture, and the other is equivalent to the union of the 
Geometrization Conjecture and Thurston's Virtual Bundle Conjecture. 
\end{abstract}

\maketitle

\section{Introduction}

The Poincar\'{e} Conjecture states that every closed, simply connected 3-manifold is 
homeomorphic to $S^3$. Stallings \cite{St2} and Jaco \cite{Ja1} have shown that 
the Poincar\'{e} Conjecture is equivalent to a purely algebraic conjecture. 
Let $S$ be a closed, orientable 
surface of genus $g$ and $F_1$ and $F_2$ free groups of rank $g$. A homomorphism 
$\varphi=\varphi_1 \times \varphi_2 : \pi_1(S) \rightarrow F_1 \times F_2$ is 
called a {\it splitting homomorphism} of genus $g$ if $\varphi_1$ and $\varphi_2$ are 
onto. It has an {\it essential factorization} through a free product if 
$\varphi=\theta\circ\psi$, where 
$\psi:\pi_1(S)\rightarrow A*B$, $\theta: A*B\rightarrow F_1 \times F_2$, and $\im \psi$ 
is not conjugate into $A$ or $B$. 

\begin{thm}[Stallings-Jaco] The Poincar\'{e} Conjecture is true if and only if 
every splitting epimorphism of genus $g > 1$ has an essential factorization. 
\hfill \qedsymbol \end{thm}

Thurston's Geometrization Conjecture \cite[Conjecture 1.1]{Th} is equivalent to 
the statement that each prime connected summand of a closed, 
connected, orientable 3-manifold either is Seifert fibered, is hyperbolic, or 
contains an incompressible torus. (See \cite{Sc3}.) In particular, 
it implies that a closed, connected, orientable 
3-manifold with finite fundamental group must be Seifert fibered. Since a closed, 
simply connected Seifert fibered space must be homeomorphic to $S^3$,  
the Poincar\'{e} Conjecture is a special case of the Geometrization Conjecture. 
The goal of this paper is to generalize the Stallings-Jaco theorem to the setting 
of the Geometrization Conjecture. 

The special case in which the fundamental group is finite is the closest analogue of 
the Stallings-Jaco theorem. Denote by $[G:H]$ the index of the subgroup $H$ of the 
group $G$. 

\begin{thm} The Geometrization Conjecture is true for closed, connected, orientable 
3-manifolds with 
finite fundamental group if and only if every splitting homomorphism $\varphi$ of genus 
$g > 2$ such that $[F_1 \times F_2:\im \varphi]<\infty$ 
has an essential factorization. \end{thm}

For the general case we have the following result. 

\begin{thm} The Geometrization Conjecture is true if and only if for every splitting 
homomorphism $\varphi$ of genus $g > 2$ either 
\begin{enumerate}
\item $\varphi$ has an essential factorization, or 
\item $\pi_1(S)/\ker \varphi_1  \ker \varphi_2$ either 
\begin{enumerate}
\item contains a $\mathbf{Z}\oplus\mathbf{Z}$ subgroup, or 
\item is isomorphic to a discrete, non-trivial, torsion-free subgroup of 
$SL(2,\mathbf{C})$. 
\end{enumerate}
\end{enumerate} \end{thm}

Thurston has also conjectured \cite[p.\ 380]{Th} that every closed, connected, 
hyperbolic 3-manifold has a 
finite sheeted covering space which is a surface bundle over $S^1$. Combining this 
``Virtual Bundle Conjecture'' with the Geometrization Conjecture gives an ``Extended 
Geometrization Conjecture'' which is also equivalent to an algebraic conjecture. 
Define a subgroup $H$ of a group $G$ to be {\it good} if it is finitely generated 
and non-trivial and $N(H)/H$ has an element of infinite order, where $N(H)$ is the 
normalizer $\{g \in G \, | \, gHg^{-1}=H\}$ of $H$ in $G$. 

\begin{thm} The Extended Geometrization Conjecture is true if and only if 
for every splitting homomorphism $\varphi$ of genus $g>1$ either 
\begin{enumerate}
\item $\varphi$ has an essential factorization, or 
\item $\pi_1(S)/\ker \varphi_1 \, \ker \varphi_2$ has a good subgroup, or 
\item $[F_1 \times F_2 : \im  \varphi]< \infty$, and $g=2$.
\end{enumerate} \end{thm}

The paper is organized as follows. Section 2 quotes some background lemmas. 
Section 3 proves some algebraic results. Sections 4, 5, and 6 prove Theorems 
1.2, 1.3, and 1.4, respectively. Section 7 gives an alternative proof of 
Thurston's observation \cite[p.\ 380]{Th} (see also Culler and Shalen 
\cite[Theorem 4.2.1]{Cu-Sh} and Gabai \cite{Ga2}) that closed, orientable, virtually 
hyperbolic 3-manifolds are homotopy hyperbolic. 

Unless the contrary is evident all manifolds under consideration are assumed 
to be connected.

\section{Heegaard splittings and splitting homomorphisms} 

Recall that a \textit{Heegaard splitting} of a closed, orientable 
3-manifold $M$ is 
a pair $(M,S)$, where $S$ is a closed, orientable surface in $M$ such that $M-S$ 
has two components, and the closures of these components are cubes with handles, 
which we denote by $V_1$ and $V_2$. The \textit{genus} of the splitting is the 
genus of $S$. The splitting is \textit{reducible} if there is a 2-sphere $\Sigma$ in 
$M$ which is in general position with respect to $S$ such that $S \cap \Sigma$ is a 
simple closed curve which does not bound a disk on $S$. Recall also that $M$ is said 
to be \textit{reducible} if it contains a 2-sphere which does not bound a 3-ball. 

\begin{lem}[Haken] Every Heegaard splitting of a reducible 3-manifold is reducible. 
\end{lem}

\begin{proof} See \cite[p.\ 84]{Ha}. See also \cite[Theorem II.7]{Ja2}. \end{proof} 

It follows that if $(M,S)$ is a genus $g$ Heegaard splitting of a reducible 3-manifold 
$M$, then either $g=1$ and $M$ is homeomorphic to $S^1 \times S^2$ or $g>1$ and $M$ can 
be expressed as a connected sum of 3-manifolds having Heegaard splittings of lower genera. 
See \cite[Lemma 3.8]{He}.

Two splitting homomorphisms $\varphi:\pi_1(S) \rightarrow F_1 \times F_2$ and 
$\varphi^{\prime}:\pi_1(S^{\prime}) \rightarrow F_1^{\prime} \times F_2^{\prime}$ are 
\textit{equivalent} if there are isomorphisms $\sigma_i:F_i \rightarrow F_i^{\prime}$ 
and $\tau:\pi_1(S) \rightarrow \pi_1(S^{\prime})$ such that $\sigma \circ \varphi=
\varphi^{\prime} \circ \tau$, where 
$\sigma=\sigma_1 \times \sigma_2$. Note that in this case $\varphi$ has an essential 
factorization if and only if $\varphi^{\prime}$ does. 

Every Heegaard splitting gives rise to a splitting homomorphism by choosing 
a basepoint on $S$, and, for $i=1$, $2$, letting $F_i=\pi_1(V_i)$ and letting $\varphi_i$ 
be the induced homomorphism on fundamental groups. A splitting homomorphism is 
\textit{realized} by a Heegaard splitting if it is equivalent to a splitting 
homomorphism of this type. 

\begin{lem}[Jaco] Every splitting homomorphism can be realized by a Heegaard splitting 
of some 3-manifold. \end{lem} 

\begin{proof} This is Theorem 5.2 of \cite{Ja1}. \end{proof} 

\begin{lem}[Stallings-Jaco] Suppose $g>1$. Then $(M,S)$ is reducible if and only 
if the associated splitting homomorphism has an essential factorization. \end{lem}

\begin{proof} Sufficiency is due to Stallings \cite[Theorem 2]{St2} and necessity 
to Jaco \cite[pp.\ 377--378]{Ja1}. \end{proof} 

\begin{lem}[Stallings] $\pi_1(M)$ is isomorphic to 
$\pi_1(S)/\ker \varphi_1 \, \ker \varphi_2$. \end{lem}

\begin{proof} See \cite[p.\ 85]{St2}. See also \cite[p.\ 128]{Ma}.  \end{proof}

We briefly sketch how these ingredients give the Stallings-Jaco theorem. 
Stallings showed that $\pi_1(M)$ is trivial if and only if $\varphi$ is onto 
\cite[Theorem 1]{St2}. (See also Lemma 3.1 below.)  
Thus if every splitting epimorphism of genus greater than one has an essential 
factorization then every Heegaard splitting of genus greater than one of a 
homotopy 3-sphere is reducible, and so one can express it as a connected sum 
of homotopy 3-spheres with genus one Heegaard splittings, which must be 
homeomorphic to $S^3$. Jaco took an arbitrary splitting epimorphism of genus 
greater than one and realized it by a Heegaard splitting of a homotopy 3-sphere. 
If it is homeomorphic to $S^3$, then by a result of Waldhausen \cite{Wa} the 
splitting is reducible, and hence the splitting epimorphism has an essential 
factorization. 

We finally remark that the Geometrization Conjecture is well known to hold for 
closed, orientable 3-manifolds with Heegaard splittings of genus at most two. 
For genus at most one this is classical. For genus two, such manifolds admit 
involutions with 1-manifolds as fixed point sets \cite{Bi-Hi}. The result then 
follows from Thurston's Orbifold Theorem, which was announced in \cite[p.\ 362]{Th}, 
and has been given detailed proofs by Cooper, Hodgson, and Kerckhoff and also by 
Boileau and Porti \cite{Boi-Po} in the case of a cyclic group action with 
1-dimensional fixed point set. 

\section{Some algebraic results}

For a subgroup $H$ of a group $G$ let $G/H$ denote the set of left cosets $[g]=gH$ of 
$H$ in $G$. 

\begin{lem} Let $\varphi$ be a splitting homomorphism. There is a bijection 
\[ \Phi:\pi_1(M)\cong\pi_1(S)/\ker \varphi_1 \, \ker \varphi_2 \rightarrow 
F_1 \times F_2 \, /\im  \varphi \] 
given by $\Phi([x])=
[(\varphi_1(x),1)]$. \end{lem} 

\begin{proof} $\Phi$ is well defined: Suppose $y=xk_1k_2$, where 
$k_i \in \ker \varphi_i$. Then 
\begin{eqnarray*} 
\Phi([y]) & = & [(\varphi_1(xk_1k_2),1)]\\ 
 & = & [(\varphi_1(x)\varphi_1(k_1)\varphi_1(k_2),1)]\\
 & = & [(\varphi_1(x)\varphi_1(k_2),1)]\\
 & = & [(\varphi_1(x)\varphi_1(k_2),\varphi_2(k_2))]\\ 
 & = & [(\varphi_1(x),1)(\varphi_1(k_2),\varphi_2(k_2))]\\
 & = & [(\varphi_1(x),1)\varphi(k_2)]\\
 & = & [(\varphi_1(x),1)]\\
 & = & \Phi([x]) \end{eqnarray*} 

$\Phi$ is one to one: If $\Phi([x])=\Phi([y])$, 
then for some $z \in \pi_1(S)$ one has 
\begin{eqnarray*} (\varphi_1(x),1) & = &(\varphi_1(y),1) \varphi(z)\\
 & = & (\varphi_1(y),1)(\varphi_1(z),\varphi_2(z))\\
 & = & (\varphi_1(y)\varphi_1(z),\varphi_2(z)) \end{eqnarray*} 
Thus $z \in \ker \varphi_2$, and $\varphi_1(x)=\varphi_1(yz)$. Let $k_1=x(yz)^{-1}$. 
Then $k_1 \in \ker \varphi_1$, and $x=k_1yz=yk_1^{\prime}z$ for some $k_1^{\prime} 
\in \ker \varphi_1$ since this subgroup is normal in $\pi_1(S)$. Thus $[x]=[y]$. 

$\Phi$ is onto: Let $(a,b) \in F_1 \times F_2$. Since the $\varphi_i$ are onto we have 
$a=\varphi_1(x)$ and $b=\varphi_2(y)$ for some $x$, $y \in \pi_1(S)$. Thus 
\begin{eqnarray*} [(a,b)] & = & [(\varphi_1(x),\varphi_2(y))]\\
 & = & [(\varphi_1(xy^{-1})\varphi_1(y),\varphi_2(y)]\\
 & = & [(\varphi_1(xy^{-1}),1)(\varphi_1(y),\varphi_2(y))]\\
 & = & [(\varphi_1(xy^{-1}),1)]\\
 & = & \Phi([xy^{-1}]) \end{eqnarray*} \end{proof}

In general $\Phi$ need not be an isomorphism because $\im \varphi$ need not be normal in 
$F_1 \times F_2$, and so $F_1 \times F_2 \, /\im  \varphi$ need not be a group. In fact, 
one has the following precise result. Let $Z(G)$ denote the center of the group $G$. 
Recall that $N(H)$ denotes the normalizer of the subgroup $H$ of $G$. 

\begin{prop} $\Phi(Z(\pi_1(M)))=N(\im \varphi)/\im \varphi$. \end{prop}

\begin{proof} Suppose $[y] \in Z(\pi_1(M))$. 
Let $x \in \pi_1(S)$. Then $yxy^{-1}x^{-1}=k_1k_2$ for some $k_i \in \ker \varphi_i$, 
$i=1$, $2$. So $\varphi_1(yxy^{-1})=\varphi_1(yxy^{-1}x^{-1}x)=\varphi_1(k_1k_2x)=
\varphi_1(k_2x)$ and $\varphi_2(x)=\varphi_2(k_2x)$. Thus 
\[(\varphi_1(y),1)(\varphi_1(x),\varphi_2(x))(\varphi_1(y^{-1}),1)=
(\varphi_1(yxy^{-1}),\varphi_2(x))=(\varphi_1(k_2x),\varphi_2(k_2x)).\] 
Similarly $y^{-1}xyx^{-1}=k_1^{\prime}k_2^{\prime}$ implies that 
\[(\varphi_1(y^{-1}),1)(\varphi_1(x),\varphi_2(x))(\varphi_1(y),1)=
(\varphi_1(k_2^{\prime}x),\varphi_2(k_2^{\prime}x)).\] 
Hence $(\varphi_1(y),1) \in N(\im \varphi)$. 

Now suppose that $(a,b) \in N(\im \varphi)$. From the proof that $\Phi$ is onto in 
Lemma 3.1 we may assume that $(a,b)=(\varphi_1(y),1)$ for some $y \in \pi_1(S)$. 
Let $x \in \pi_1(S)$. Then $(\varphi_1(y),1)(\varphi_1(x),\varphi_2(x))
(\varphi_1(y^{-1}),1)=(\varphi_1(z),\varphi_2(z))$ for some $z \in \pi_1(S)$. 
So $\varphi_1(yxy^{-1})=\varphi_1(z)$ and $\varphi_2(x)=\varphi_2(z)$. 
Hence $yxy^{-1}z^{-1}=k_1 \in \ker \varphi_1$ and $zx^{-1}=k_2 \in \ker \varphi_2$. 
So $yxy^{-1}x^{-1}=k_1k_2 \in \ker \varphi_1 \ker \varphi_2$. Hence 
$[y] \in Z(\pi_1(M))$. \end{proof}

\begin{cor} $\im  \varphi$ is normal in $F_1 \times F_2$ if and only if $\pi_1(M)$ 
is abelian.  \hfill \qedsymbol \end{cor} 

\section{The finite fundamental group case}

\begin{proof}[Proof of Theorem 1.2] First assume that every splitting homomorphism 
$\varphi$ with $g > 2$ and $[F_1 \times F_2 : \im \varphi] < \infty$ has an essential 
factorization. Let $M$ be a closed, orientable 3-manifold with $\pi_1(M)$ finite. 
We may assume that $M$ is irreducible. Let $(M,S)$ be a Heegaard splitting of $M$ 
of minimal genus $g$. Let $\varphi$ be the associated splitting homomorphism. 
By Lemma 3.1 $[F_1 \times F_2 : \im \varphi] < \infty$, and so if $g$ were greater than 
two, then $\varphi$ would have an essential factorization, and hence $(M,S)$ would be 
reducible. Since $M$ is irreducible this would yield a Heegaard splitting of lower 
genus, contradicting the choice of $g$. Thus $g \leq 2$, and we are done. 

Now assume that the Geometrization Conjecture holds in the finite fundamental group case. 
Let $\varphi$ be a splitting homomorphism with $[F_1 \times F_2 : \im \varphi] < \infty$ 
and genus $g > 2$. Realize $\varphi$ by a Heegaard splitting $(M,S)$. Then by Lemma 3.1 
$\pi_1(M)$ is finite. Suppose $\varphi$ does not have an essential factorization. 
Then $(M,S)$ is irreducible. 

By the Geometrization Conjecture $M$ is a Seifert fibered space. Since $\pi_1(M)$ is 
finite $M$ has a Seifert fibration over a 2-sphere with at most three exceptional fibers.  
(See \cite[Theorem 12.2]{He} or \cite[p.\ 92]{Ja2}. Note that it may have a Seifert 
fibration over a projective plane with one exceptional fiber, but then it also has 
a Seifert fibration of the given type.)  
If there were fewer than three exceptional fibers, then $M$ would be $S^3$ or a lens 
space. But by results of Waldhausen \cite{Wa} and of Bonahon and Otal \cite{Bon-Ot} 
the irreducible Heegaard splittings of these spaces have, respectively, 
genus zero and one, contradicting our choice of $g$. Thus there 
are three exceptional fibers $f_i$ of multiplicities $\alpha_i >1$, $i=1$, $2$, $3$. 
Moreover, up to ordering, 
$(\alpha_1, \alpha_2, \alpha_3)$ must be one of $(2, 2, \alpha_3)$, $\alpha_3 \geq 2$, 
$(2, 3, 3)$, $(2, 3, 4)$, or $(2, 3, 5)$. 

We now recall two constructions for Heegaard splittings of closed, orientable Seifert 
fibered spaces over orientable base surfaces. 
For simplicity we restrict to the special case at hand. See \cite{Mori-Sch} and \cite{Se} 
for the general case and a more detailed description. 

First choose two of the three 
exceptional fibers. Join their image points in the base 2-sphere by an arc which misses 
the image point of the other exceptional fiber. Lift this arc to an arc in $M$ joining 
the two chosen exceptional fibers. A regular neighborhood $V$ of the resulting graph is 
a cube with two handles. It turns out that the closure of its complement is also a cube 
with handles, and so $(M,\bd V)$ is a genus two Heegaard splitting of $M$. It is called 
a \textit{vertical} Heegaard splitting. We remark that in the general case all vertical 
Heegaard splittings have the same genus $g_v$. 

Next choose one exceptional fiber $f_i$, and let $N$ be a regular neighborhood of it. The 
closure $M_0$ of the complement of $N$ is bundle over $S^1$ with fiber a surface $F$ 
\cite[Theorem 12.7]{He}, \cite[Theorem VI.32]{Ja2}. Moreover, $F$ is a branched covering 
space of the base surface of the Seifert fibration; the branch points are the images of 
the exceptional fibers and have branching indices equal to the indices of the exceptional 
fibers. Suppose $\bd F$ is connected and has intersection number $\pm 1$ with a 
meridian of the solid torus $N$. Let $H$ be a regular neighborhood of $F$ in $M_0$. 
Then $H$ is a cube with handles whose genus is twice the genus of $F$. 
It turns out that the closure of the complement of $H$ in $M$ is also a cube with 
handles, and thus $(M,\bd H)$ is a Heegaard splitting of $M$. It is called a 
\textit{horizontal} Heegaard splitting at $f_i$. Denote its genus by $g_h(f_i)$. Note 
that if either of the two conditions on $\bd F$ is violated, then by definition $M$ does 
not have a horizontal Heegaard splitting at $f_i$. Let $d$ be the least common multiple 
of $\alpha_j$ and $\alpha_k$, where $f_j$ and $f_k$ are the other two exceptional fibers. 

Moriah and Schultens \cite{Mori-Sch} have shown that every irreducible Heegaard splitting 
of a closed, orientable Seifert fibered space over an orientable base surface is either 
vertical or horizontal. Since $g>2$ and $g_v=2$ our splitting $(M,S)$ must be 
horizontal. 

Sedgwick \cite{Se} has determined precisely which vertical and horizontal Heegaard 
splittings are irreducible. In particular a horizontal splitting is irreducible if 
and only if either $g_h(f_i) \leq g_v$ or $\alpha_i > d$. In our case the first 
condition is impossible, and the second condition holds only for $(2, 2, \alpha_3)$, 
where $\alpha_3 > 2$. But in this case $M_0$ is Seifert fibered over a disk with 
two exceptional fibers of index two and so must be a twisted $I$-bundle over a 
Klein bottle; it follows that the fiber $F$ is an annulus, and so $M$ does not have 
a horizontal Heegaard splitting at $f_3$.  

Thus $(M,S)$ must be reducible, and so $\varphi$ must have an essential 
factorization. \end{proof}

\section{The general case}

\begin{proof}[Proof of Theorem 1.3] Suppose the condition on the $\varphi$ holds. Let $M$ 
be a closed, orientable, irreducible 3-manifold and $(M,S)$ a Heegaard splitting of 
minimal genus $g$. If $g \leq 2$, then the Geometrization Conjecture holds for $M$. 
So assume $g>2$. If $\varphi$ had an essential factorization, then $(M,S)$ would be 
reducible, and, since $M$ is irreducible $M$ would have a Heegaard splitting of lower 
genus, contradicting the choice of $g$. So $\varphi$ does not have an essential 
factorization, and we must be in case (2). 

In case (2)(a) $\pi_1(M)$ has a $\mathbf{Z}\oplus\mathbf{Z}$ subgroup. By Scott's 
version of the torus theorem \cite{Sc2} either $M$ contains an incompressible torus,  
and we are done, or $\pi_1(M)$ contains a normal $\mathbf{Z}$ subgroup. In the latter 
case the proof of the Seifert fibered space conjecture by Casson and Jungreis 
\cite{Ca-Ju} and by Gabai \cite{Ga1} gives that $M$ is a Seifert fibered space, and 
again we are done. 

In case (2)(b) $\pi_1(M)$ is isomorphic to a discrete, non-trivial, torsion-free 
subgroup of $SL(2,\mathbf{C})$. Since the kernel of the projection to 
$PSL(2,\mathbf{C})=SL(2,\mathbf{C})/\{\pm I\}$ is a finite group this subgroup 
projects isomorphically to a discrete subgroup $\Gamma$ of $PSL(2,\mathbf{C})$. 
A subgroup of $PSL(2,\mathbf{C})$ is discrete and torsion free if and only if 
its natural action on hyperbolic 3-space 
$\mathbf{H}^3$ is free (no non-trivial element has a fixed point) and discontinuous 
(each compact set meets only finitely many of its translates) 
\cite[Theorems 8.2.1, 8.1.2, and 5.3.5]{Ra}. 
Thus the quotient space $N=\mathbf{H}^3/\Gamma$ is a hyperbolic 3-manifold 
\cite[Theorem 8.1.3]{Ra}. Since $M$ is closed and aspherical 
$H_3(N)\cong H_3(M) \cong \mathbf{Z}$, and so $N$ is closed and orientable. By the 
topological rigidity of hyperbolic 3-manifolds, due to Gabai, Meyerhoff, and N. Thurston 
\cite{Ga-Me-Th} we have that $M$ and $N$ are homeomorphic, and we are done. 

Now suppose that the Geometrization Conjecture is true. Let $\varphi$ be a splitting 
homomorphism of genus $g>2$. Assume that $\varphi$ has no essential factorization. 
Let $(M,S)$ be a Heegaard splitting which realizes $\varphi$. Then $(M,S)$ is 
irreducible, and hence $M$ is irreducible. By Lemma 3.1 and Theorem 1.2 we may assume 
that $\pi_1(M)$ is infinite. Since $M$ is irreducible, $\pi_1(M)$ is torsion-free 
\cite[Corollary 9.9]{He}.

If $M$ is hyperbolic, then $M=\mathbf{H}^3/\Gamma$ for some subgroup $\Gamma$ 
of $PSL(2,\mathbf{C})$ acting freely and discontinuously on $\mathbf{H}^3$. 
Thus $\Gamma$ is discrete and torsion free. 
By a result of Thurston $\Gamma$ lifts isomorphically to 
a subgroup of $SL(2,\mathbf{C})$. (See Culler and Shalen 
\cite[Proposition 3.1.1]{Cu-Sh}.) 

If $M$ contains an incompressible torus, then clearly $\pi_1(M)$ contains a 
$\mathbf{Z} \oplus \mathbf{Z}$ subgroup. If $M$ is Seifert fibered, then $\pi_1(M)$ 
infinite implies that $M$ has a covering space which is homeomorphic to an 
$S^1$ bundle over a closed, orientable surface $F$ of positive genus \cite[p.\ 438]{Sc3}, 
and so again $\pi_1(M)$ contains a $\mathbf{Z} \oplus \mathbf{Z}$ subgroup. 
(See also \cite[p.\ 477]{Sc3}.) \end{proof}

\section{The extended conjecture}

Recall that a non-trivial subgroup $H$ of a group $G$ is \textit{good} if it is 
finitely generated and $N(H)/H$ has an element of infinite order. 

\begin{lem} Let $M$ be a closed, orientable, irreducible 3-manifold. Then $\pi_1(M)$ 
has a good subgroup if and only if either $M$ has a finite sheeted covering 
space which is a surface bundle over $S^1$ or $\pi_1(M)$ contains a $\mathbf{Z} \oplus 
\mathbf{Z}$ subgroup. \end{lem}

\begin{proof} If $M$ is covered by a bundle with fiber a surface $F$, then the image of 
$\pi_1(F)$ in $\pi_1(M)$ is a good subgroup. If $\pi_1(M)$ contains a $\mathbf{Z} \oplus 
\mathbf{Z}$ subgroup, then a summand is a good subgroup. 

The converse follows from \cite[Lemma 1]{My}. For convenience we give the relevant 
portion of the proof of that result.  
Let $H$ be a good subgroup and $\widetilde{M}$ the covering of $M$ corresponding to $H$. 
Since the group of covering translations is isomorphic to $N(H)/H$ 
\cite[Corollary 7.3]{Ma} there is 
a covering translation of infinite order. Let $M^*$ be the quotient of 
$\widetilde{M}$ by this covering translation. $\pi_1(M^*)$ has a normal subgroup which is 
isomorphic to $H$ and has infinite cyclic quotient. The Scott compact core \cite{Sc1} $C$ 
of $M^*$ is a compact submanifold of $M^*$ with $\pi_1(C)$ isomorphic to $\pi_1(M^*)$. 
Since $M^*$ is irreducible \cite{MSY} we may assume that $C$ is irreducible. It then 
follows from the Stallings fibration theorem \cite{St1} that $C$ is a surface bundle over 
$S^1$. If $C=M^*$, then we are done. If $C \neq M^*$, then $\bd C$ consists of tori 
which are incompressible in $M^*$ and so $\pi_1(M)$ has a $\mathbf{Z} \oplus \mathbf{Z}$ 
subgroup. \end{proof}

\begin{proof}[Proof of Theorem 1.4] Suppose the extended conjecture is true. If $g>2$, 
then by the proof of Theorem 1.3 we reduce to the situation in which either $\pi_1(M)$ 
has a $\mathbf{Z} \oplus \mathbf{Z}$ subgroup or $M$ is hyperbolic; in the latter case  
we apply the Virtual Bundle Conjecture and so conclude in both cases that $\pi_1(M)$ has 
a good subgroup. If $g=2$, then a similar argument shows that either $\pi_1(M)$ has a 
good subgroup or $M$ is a Seifert fibered space with $\pi_1(M)$ finite, in which case 
$[F_1 \times F_2 : \im \varphi]<\infty$.

Now assume that the conditions on the $\varphi$ hold. 
Let $M$ be a closed, orientable, irreducible 3-manifold. By arguments similar to those 
in the proof of Theorem 1.3 we reduce to the case that $\pi_1(M)$ has a good subgroup and 
does not have a $\mathbf{Z} \oplus \mathbf{Z}$ subgroup. Then by 
Lemma 6.1 $M$ is finitely covered by a surface bundle. Note that this is the case if $M$ 
is assumed to be hyperbolic \cite[p.\ 52]{Mor}, \cite[Corollary 4.6]{Sc3}, and so the 
Virtual Bundle Conjecture holds. 
For the general situation note that by the fibered case of Thurston's hyperbolization 
theorem \cite[Theorem 2.3]{Th}, \cite{Ot} the surface bundle is hyperbolic. We conclude 
from the following result that the Geometrization Conjecture holds for $M$.

\begin{lem} Let $M$ be a closed, orientable 3-manifold which has a finite sheeted 
covering space $M^*$ which is hyperbolic. Then $M$ is hyperbolic. \end{lem} 

\begin{proof} This is well known. It follows immediately from the topological rigidity 
of hyperbolic 3-manifolds \cite{Ga-Me-Th} and the observation of Thurston that $M$ is 
homotopy equivalent to a hyperbolic 3-manifold. See \cite[Theorem 4.2.1]{Cu-Sh} or the 
next section for a proof. \end{proof}

This concludes the proof of Theorem 1.4. \end{proof}

\section{Virtually hyperbolic 3-manifolds}

The following result was observed by Thurston 
\cite[p.\ 380]{Th} to be a consequence of the Mostow rigidity theorem \cite{Mos}. 
A proof was given by Culler and Shalen \cite[Theorem 4.2.1]{Cu-Sh}; a sketch of 
the proof has also been given by Gabai \cite{Ga2}. In this section we fill in 
some details of this sketch to give a proof which, though similar to that of 
Culler and Shalen, makes somewhat less explicit use of hyperbolic geometry. 

\begin{lem}[Thurston] Let $M$ be a closed, orientable 3-manifold which 
has a finite sheeted covering space which is hyperbolic. Then $M$ is homotopy 
equivalent to a hyperbolic 3-manifold. \end{lem}

\begin{proof}  
We may assume that the covering is regular 
\cite[Theorem 4.7]{Ly-Schu}. 

As pointed out by Gabai, the covering translation in $M^*$ corresponding to an element 
of $\pi_1(M)$ is by Mostow rigidity homotopic to a unique isometry. The lifts of these 
isometries to the universal cover $\mathbf{H}^3$ give a subgroup $\Gamma$ of 
$I\!som(\mathbf{H}^3)$. The quotient $N=\mathbf{H}^3/\Gamma$ is then  
the desired hyperbolic 3-manifold. To fill out this sketch one needs to verify that 
$\pi_1(M)\cong\Gamma$ and that $\Gamma$ acts freely and discontinuously on 
$\mathbf{H}^3$. It will be convenient to first establish the following result. 

\begin{lem} Let $M^*$ be a closed hyperbolic 3-manifold and $H:M^* \times I 
\rightarrow M^*$ 
a homotopy such that $H(x,0)=H(x,1)=x$ for all $x \in M^*$. Fix $y_0^* \in M^*$. Let 
$m(t)=H(y_0^*,t)$ for all $t \in I$. Then the class $\mu$ of $m$ in $\pi_1(M^*,y_0^*)$ is 
trivial. \end{lem}

\begin{proof} Let $\lambda \in \pi_1(M^*,y_0^*)$ be represented by the loop $\ell(s)$. 
Then the map $G:I \times I \rightarrow M^*$ given by $G(s,t)=H(\ell(s),t)$ shows that 
$\mu \lambda=\lambda \mu$. Hence $\mu \in Z(\pi_1(M^*,y_0^*))$, which is trivial for a 
closed hyperbolic 3-manifold \cite[p.\ 52]{Mor}. \end{proof}

Returning to the proof of Lemma 7.1, we have covering maps 
$\mathbf{H}^3 \stackrel{p}{\rightarrow} M^* \stackrel{q}{\rightarrow} M$. 
Choose a basepoint $\tilde{x}_0 \in \mathbf{H}^3$, and let $x_0^*=p(\tilde{x}_0)$ and 
$x_0=q(x_0^*)$. Given $\alpha \in \pi_1(M,x_0)$, let $\alpha^*$ and $\tilde{\alpha}$ be 
path classes lifting $\alpha$ with $\alpha^*(0)=x_0^*$ and $\tilde{\alpha}(0)=\tilde{x}_0$, 
respectively. There is a covering translation $f_0$ of $q$ such that 
$f_0(x_0^*)=\alpha^*(1)$ and a lifting $\tilde{f}_0$ of $f_0$ such that 
$\tilde{f}_0(\tilde{x}_0)=\tilde{\alpha}(1)$. By Mostow rigidity there is a unique 
isometry $f_1$ of $M^*$ which is homotopic to $f_0$. Let $f_t$ be a homotopy from 
$f_0$ to $f_1$. It lifts to a homotopy $\tilde{f}_t$ of $\tilde{f}_0$ to an isometry 
$\tilde{f}_1$ of $\mathbf{H}^3$. 

We claim that $\tilde{f}_1$ is independent of the choice of homotopy $f_t$. 
Let $f^{\prime}_t$ be another homotopy from $f_0$ to $f_1$. Define $H: M^* \times I 
\rightarrow M^*$ by $H(x,t)=f_{2t}(f_0^{-1}(x))$ for $t \in [0,\frac{1}{2}]$ and 
$H(x,t)=f_{2-2t}^{\prime}(f_0^{-1}(x))$ for $t \in [\frac{1}{2},1]$. This homotopy 
satisfies the conditions of Lemma 7.2, and so the loop $H(f_0(x_0^*),t)$ is homotopically 
trivial. This implies that the paths $f_t(x_0^*)$ and $f_t^{\prime}(x_0^*)$ are 
path homotopic, and so their liftings $\tilde{f}_t(\tilde{x}_0)$ and 
$\tilde{f}_t(\tilde{x}_0)$ have the same endpoint $\tilde{f}_1(\tilde{x})=
\tilde{f}_1^{\prime}(\tilde{x})$. Thus $\tilde{f}_1=\tilde{f}_1^{\prime}$. 

Thus the function $\Psi:\pi_1(M,x_0) \rightarrow I\!som(\mathbf{H}^3)$ given by 
$\Psi(\alpha)=\tilde{f}_1$ is well defined. Note that if $\alpha \in \im q_*$, 
then $f_0=f_1=id_{M^*}$, and $\tilde{f}_0=\tilde{f}_1$. Let $\Gamma=\im \Psi$ and 
$\Gamma_0=\Psi(\im q_*)$. 

We next show that $\Psi$ is a homomorphism. 
Suppose $\Psi(\beta)=\tilde{g}_1$ and $\Psi(\alpha \beta)=\tilde{h}_1$. 
Let $\gamma^*$ be the image of $\beta^*$ under $f_0$. Then $\alpha^* \gamma^*=
(\alpha \beta)^*$, and so $f_0(g_0(x_0^*))=f_0(\beta^*(1))=\gamma^*(1)=h_0(x_0^*)$. 
Thus $f_0 \circ g_0=h_0$, and $f_t \circ g_t$ is a homotopy from this map to the 
isometry $f_1 \circ g_1$. By the uniqueness of $h_1$ we have that $f_1 \circ g_1=h_1$. 
By choosing the homotopy $h_t$ from $h_0$ to $h_1$ to be $f_t \circ g_t$, we see that 
$\tilde{h}_1=\tilde{f}_1 \circ \tilde{g}_1$, and so $\Psi(\alpha \beta)=\Psi(\alpha) 
\Psi(\beta)$. 

We now show that $\Psi$ is one to one. Suppose $\Psi(\alpha)=\tilde{f}_1=
id_{\mathbf{H}^3}$. If $\alpha \in \im q_*$, then $\tilde{f}_0=\tilde{f}_1$, and 
so $\alpha$ is trivial. So assume that this is not the case. Then $f_1$ is the identity 
of $M^*$, and $f_0$ is not. Since the covering is finite sheeted $f_0^n$ is the 
identity for some $n>0$. Define $H: M^* \times I \rightarrow  M^*$ by 
$H(x,t)=f_0^{k-1}(f_{k-nt}(x))$ for $t \in [\frac{k-1}{n}, \frac{k}{n}]$. Then 
$H(x,\frac{k}{n})=f_0^k$. By Lemma 7.2 the class $\mu$ of the loop $m(t)=H(x_0^*,t)$ is 
trivial in $\pi_1(M^*, x_0^*)$. Let $\rho$ be the path class of the path 
$r(t)=f_{1-t}(x_0^*)$ joining $x_0^*$ and $f_0(x_0^*)$. Then $q_*(\rho)$ is non-trivial, 
but $(q_*(\rho))^n=q_*(\mu)$ is trivial. Thus $\pi_1(M,x_0)$ has torsion, 
contradicting the fact that $M$ is aspherical \cite[Corollary 9.9]{He}. 

$\Gamma$ is discrete because it contains the finite index discrete subgroup 
$\Gamma_0$ \cite[Lemma 8, p. 177]{Ra}. Thus it acts discontinuously on $\mathbf{H}^3$. 
It acts freely because otherwise it would have torsion, contradicting the asphericity 
of $M$. Hence $N=\mathbf{H}^3/\Gamma$ is a 3-manifold with $\pi_1(N)\cong\pi_1(M)$. 
It follows from asphericity that $N$ and $M$ are homotopy equivalent, and so $N$ is 
closed and orientable; hence $\Gamma \subset PSL(2,\mathbf{C})$. \end{proof}


\begin{thebibliography}{99}

\bibitem{Bi-Hi}
J. S. Birman and H. M. Hilden, {\it Heegaard splittings of branched coverings of 
$S^3$}, Trans. Amer. Math. Soc. 213 (1975), 315--352. 

\bibitem{Boi-Po}
M. Boileau and J. Porti, {\it Geometrization of 3-orbifolds of cyclic type}, 
xxx archive preprint math.GT/9805073.

\bibitem{Bon-Ot}
F. Bonahon and J-P. Otal, {\it Scindements de Heegaard des espaces lenticulaires}, 
Ann. Sci. \'{E}cole Norm. Sup. (4) 16 (1983), no. 3, 451--466 (1984).

\bibitem{Ca-Ju}
A. Casson and D. Jungreis, {\it Convergence groups and Seifert fibered 3-manifolds},
Invent. Math. 118 (1994), 441--456.

\bibitem{Cu-Sh}
M. Culler and P. B. Shalen, {\it Varieties of group representations and splittings of 
$3$-manifolds}, Ann. of Math. (2) 117 (1983), no. 1, 109--146.

\bibitem{Ga1}
D. Gabai, {\it Convergence groups are Fuchsian groups},
Annals of Math.  136 (1992), 447--510.

\bibitem{Ga2}
D. Gabai, {\it Homotopy hyperbolic 3-manifolds are virtually hyperbolic}, 
J. Amer. Math. Soc. 7 (1994), 193--198 

\bibitem{Ga-Me-Th} 
D. Gabai, R. Meyerhoff, and N. Thurston, {\it Homotopy hyperbolic 
3-manifolds are hyperbolic}, MSRI Preprint Series \#1996-058. 

\bibitem{Ha}
W. Haken, {\it Some results on surfaces in $3$-manifolds}, 1968 Studies in Modern 
Topology pp. 39--98, Math. Assoc. Amer. (distributed by Prentice-Hall, Englewood 
Cliffs, N.J.)

\bibitem{He}
J. Hempel, {\it 3-Manifolds}, Ann. of Math. Studies, No. 86, Princeton,
(1976).

\bibitem{Ja1}
W. Jaco, {\it Heegaard splittings and splitting homomorphisms}, 
Trans. Amer. Math. Soc. 144 (1969), 365--379.

\bibitem{Ja2}
W. Jaco, {\it Lectures on three-manifold topology}, CBMS Regional
Conference Series in Math., No. 43, Amer. Math. Soc. (1980).

\bibitem{Ly-Schu}
R. Lyndon and P. Schupp, Combinatorial Group Theory, Ergebnisse der Mathematik und ihrer
Grenzgebiete, Band 89. Springer-Verlag, Berlin-New York, 1977.

\bibitem{Ma}
W. Massey, {\it Algebraic Topology: An Introduction}, Graduate Texts in 
Mathematics No. 56, Springer-Verlag, 1977. 

\bibitem{MSY}
W. Meeks, L. Simon, S. T. Yau, {\it Embedded minimal surfaces, exotic 
spheres, and manifolds with postive Ricci curvature}, Annals of Math, 
116 (1982), 621--659.

\bibitem{Mor}
J. W. Morgan, {\it On Thurston's uniformization theorem for three-dimensional manifolds}, 
The Smith Conjecture (New York, 1979), 37--125, Pure Appl. Math., 112, Academic Press, 
Orlando, FL, 1984.

\bibitem{Mori-Sch}
Y. Moriah and J. Schultens, {\it Irreducible Heegaard splittings of Seifert fibered 
spaces are either vertical or horizontal}, Topology 37 (1998), no. 5, 1089--1112.

\bibitem{Mos}
G. D. Mostow, {\it Quasi-conformal mappings in $n$-space and the rigidity of hyperbolic 
space forms}, Inst. Hautes Études Sci. Publ. Math. No. 34 (1968), 53--104.

\bibitem{My}
R. Myers, {\it Compactifying sufficiently regular covering spaces of compact 3-manifolds}, 
Proc. Amer. Math. Soc., to appear. 

\bibitem{Ot}
J-P. Otal, {\it Le th\'{e}or\`{e}me d'hyperbolisation pour les 
vari\'{e}t\'{e}s fibr\'{e}es 
de dimension 3,} Ast\'{e}risque No. 235 (1996),

\bibitem{Ra}
J. Ratcliffe, {\it Foundations of hyperbolic manifolds}, 
Graduate Texts in Mathematics, 149. Springer-Verlag, New York, 1994.

\bibitem{Sc1} 
P. Scott, {\it Compact submanifolds of 3-manifolds}, J. London Math. 
Soc. 7 (1973), 246--250. 

\bibitem{Sc2}
P. Scott, {\it A new proof of the annulus and torus theorems}, 
Amer. J. Math. 102 (1980), no. 2, 241--277. 

\bibitem{Sc3}
P. Scott, {\it The geometries of $3$-manifolds}, 
Bull. London Math. Soc. 15 (1983), no. 5, 401--487.

\bibitem{Se}
E. Sedgwick, {\it The irreducibility of Heegaard splittings of Seifert 
fibered spaces}, Pacific J. Math., to appear. 

\bibitem{St1}
J. Stallings, {\it On fibering certain $3$-manifolds}, 1962 Topology 
of 3-manifolds and related topics (Proc. The Univ. of Georgia Institute,
1961) pp. 95--100 Prentice-Hall, Englewood Cliffs, N.J.

\bibitem{St2}
J. Stallings, {\it How not to prove the Poincar\'{e} conjecture}, Topology 
Seminar, Wisconsin, 1965, Ann. of Math. Studies, No. 60, Princeton,
(1966).

\bibitem{Th}
W. P. Thurston, {\it Three-dimensional manifolds, Kleinian groups and 
hyperbolic geometry}, Bull. Amer. Math. Soc. (N.S.) 6 (1982), no. 3,
357--381.

\bibitem{Wa}
F. Waldhausen, {\it Heegaard-Zerlegungen der $3$-Sph\"{a}re}, Topology 7 
(1968), 195--203.

\end{thebibliography}
\end{document}